\documentclass[submission,pdftex,10pt]{dmtcs}

\usepackage{amsmath}
\usepackage{amssymb}
\usepackage{times}
\usepackage{hyperref}
\usepackage{graphicx}

\newtheorem{proposition}{Proposition}

\newtheorem{theorem}{Theorem}
\newtheorem{lemma}{Lemma}

\def\R{\mathbb{R}}

\def\Pr{\mathbb{P}}
\def\cal{\mathcal}
\def\wh{\widehat}

\def\cau{\frac{1}{2i\pi}}

\author{Nicolas Broutin and Philippe Flajolet}
\date{\today}
\title[Height of random binary unlabelled trees]{The height of random binary unlabelled trees}
\address{Algorithms Project, INRIA-Rocquencourt, F-78153 Le Chesnay (France)}

\usepackage[small]{caption2}
\usepackage{booktabs}
\usepackage[sort,numbers]{natbib}
\usepackage{paralist}

\hypersetup{
    bookmarks=true,         
    unicode=false,          
    pdftoolbar=true,        
    pdfmenubar=true,        
    pdffitwindow=true,      
    pdftitle={My title},    
    pdfauthor={Author},     
    pdfsubject={Subject},   
    pdfnewwindow=true,      
    pdfkeywords={keywords}, 
    colorlinks=true,       
    linkcolor=blue,          
    citecolor=blue,        
    filecolor=blue,      
    urlcolor=blue           
}

\newcommand\p[1]{\mathbb{P}\left\{#1\right\}}
\newcommand{\pran}[1]{{\left(#1\right)}}

\newcommand{\E}[1]{{\mathbb E}\left[ #1\right]}

\makeatletter
\def\timenow{\@tempcnta\time
  \@tempcntb\@tempcnta
  \divide\@tempcntb60
  \ifnum10>\@tempcntb0\fi\number\@tempcntb
  \multiply\@tempcntb60
  \advance\@tempcnta-\@tempcntb:\ifnum10>\@tempcnta0\fi\number\@tempcnta}
\makeatother

\newenvironment{Proof}{\begin{proof}}{\end{proof}}

\begin{document}

\maketitle

\begin{abstract}
This extended abstract is dedicated  to the analysis  of the  height of non-plane
unlabelled  rooted  binary trees.   The height of   such a tree chosen
uniformly among those  of size~$n$ is  proved to have a limiting theta
distribution, both in a central and  local sense.  Moderate as well as
large deviations estimates  are also derived.  The proofs  rely on the
analysis  (in the complex  plane)  of generating functions  associated
with trees of bounded height.
\end{abstract}

\section{Introduction}
\renewcommand{\baselinestretch}{1}
We consider trees that are \emph{binary,  non-plane, unlabelled, and rooted};
that is, a tree is taken in the graph-theoretic sense and it has nodes
of (out)degree two or zero only; a special node is distinguished, the
root,  which   has  degree two.     In this  model,   the   nodes  are
indistinguishable, while no order is assumed between the neighbours of
a node.  Let  $\mathcal Y$ denote  the class  of  such trees, and  let
$\mathcal  Y_n$ be the  subset  consisting of  trees with $n$ external
nodes (i.e., nodes of degree~zero).

In this extended abstract, we study the (random) \emph{height} $H_n$  of a tree sampled
uniformly from $\mathcal Y_n$.  The depth of nodes 
has been analysed for many ``\emph{simple varieties}'' of trees by
Meir and Moon~\cite{MeMo78}.  Regarding height, a  few  special cases were
studied early: R\'enyi and Szekeres~\cite{ReSz67} proved 
that the average  height of \emph{labelled} non-plane  trees  of size~$n$ is
asymptotic  to  $2\sqrt{\pi n}$;  
De Bruijn, Knuth, and Rice~\cite*{BrKnRi72} dealt with
\emph{plane}  trees and  showed  that  the average height is equivalent to
$\sqrt{\pi  n}$ as   $n\to\infty$.  Finally,  Flajolet and Odlyzko~\cite{FlOd82}
developed an approach that  encompasses all simple varieties  of trees
(see~\cite{FlGaOdRi93} for more results).  These  results  are relative to trees
where one can distinguish the  neighbours of a  node, either by  their
labels (labelled trees), or by the order induced on the progeny (plane
trees).

In such models with distinguishable progeny,  there are natural random
walks associated to the random trees. Also,  trees of a fixed size $n$
may be seen as Galton--Watson  processes conditioned on the size of the total progeny being
$n$;    cf~\cite{Aldous90,Kennedy75,Kolchin86}.  
Probabilistic techniques have  been  applied to the random
walk associated to a tree traversal in order to derive   asymptotic results about   the trees,  in  particular
regarding   height and   width~\cite{ChMa01,ChMaYo00}.   Yet  an other
approach is  to find the  continuous limit of suitably rescaled random
trees of increasing sizes. One can  then read off some of the limit parameters
directly on the limiting object.  This point of  view has been adopted
by Aldous~\cite{Aldous91b}  in his definition of the continuum
random tree.  For a  recent account of the probabilistic developments,
see the survey by Le Gall~\cite{LeGall05}.

The  case   of trees  with  indistinguishable  progeny  is essentially
different, and  no direct  random  walk  approach  would  seem  to be
possible, due to the symmetries inherent in unlabelled structures. The
analysis of unlabelled non-plane trees finds its  origins in the works
of  P\'olya~\cite{Polya37} and   Otter~\cite{Otter48}.  However, these
authors mostly  focused on enumeration---the problem of characterizing
typical parameters of these   random trees remains  largely untouched.
Gittenberger~\cite{Gittenberger05} and recently, in  an  independent study that
predates ours by a couple of weeks, Drmota--Gittenberger~\cite{DrGi08} have
examined the profile  of non-plane unlabelled \emph{``general''  trees}  (where all   degrees  are
allowed) and shown that the
joint distribution of the number of nodes
at any finite number  of levels   converges  weakly  to  the  finite-dimensional   distribution of  Brownian excursion  local  times. They further extended the 
result to a convergence of the entire profile to the Brownian excursion local time. This gives in particular the limit law for the height of these trees.
This
suggests  that, although  there  is no clear  reduction of  unlabelled
trees  to random walks, such trees  largely behave like simply generated
families. In particular, this suggests that the rescaled height $H_n/\sqrt
n$    should   admit  a       limit   distribution  of    the    theta
type~\cite{FlOd82,Kennedy76,ReSz67}. We shall prove that 
such is indeed the case for \emph{binary non-plane trees} 
in Theorems~\ref{clt} and~\ref{llt} below. We also provide moderate and large deviations estimates (Theorems~\ref{thm:moderate} and~\ref{thm:large_deviations}) as well as asymptotic estimates for the moments (Theorem~\ref{thm:moments}); see~\S\ref{asy-sec}. Some of the more technical proofs are omitted and we limit ourselves to the global structure of the  arguments; the details may be found in the long version \cite{BrFl08b}.

This note  arose   from questions  of Jean-Fran\c{c}ois Marckert and Gr\'egory Miermont \cite{marckert2008}.  Their
motivation  comes from   an   attempt at extending   the probabilistic
methods of Aldous to   non-plane  trees and developing   corresponding
continuous models---we are indebted to them for being  at the origin of
the present  study.  We also express our  gratitude to  Alexis Darrasse and
Carine  Pivoteau  for  designing   and programming  for  us  efficient
Boltzmann samplers of  binary trees and providing detailed statistical
data that guided our first analyses of this problem.

\section{Trees and generating functions}\label{sec-basics}

\emph{Tree enumeration.}
Our approach is entirely based on \emph{generating functions}. The class~$\cal Y$ of
binary (non-plane unlabelled rooted) trees is defined to include the 
tree with a single external node. A tree has size~$n$ if it has~$n$ 
external nodes (hence $n-1$ internal nodes). The cardinality of 
the subclass~$\cal Y_n$ of trees of size~$n$ is denoted by~$y_n$
and the generating function (GF) of~$\cal Y$ is
\[
y(z):=\sum_{n\ge1} y_z z^n=z+z^2+z^3+2z^4+3z^5+6z^6+11z^7+23z^8+\cdots,
\]
the coefficients corresponding to the entry A001190 (Wedderburn--Etherington numbers) of Sloane's 
\emph{On-line Encyclopedia of Integer Sequences}.

Since a binary tree is either an external node or a root appended to an unordered pair of
two (not necessarily distinct) binary trees, one has the basic functional equation
\begin{equation}\label{eq:gen_eq}
y(z)=z+\frac12 y(z)^2+\frac12y(z^2),
\end{equation}
as follows from fundamental principles of combinatorial enumeration~\cite{FlSe08,HaPa73,Polya37,Otter48}.
According to the general principles of analytic combinatorics, we shall 
operate in an essential manner with properties of generating functions
in the \emph{complex plane}. The following lemma is classical:

\begin{lemma}\label{lem:rho1}
Let $\rho$ be the radius of convergence of $y(z)$. Then, 
one has $1/4 \le \rho \le 1/2$,
and~$\rho$ is determined implicitly by 
$2\rho+y(\rho^2)=1$. 
As $z\to\rho^-$, the generating function $y(z)$ satisfies
\begin{equation}\label{singy}
y(z)=1-\lambda\sqrt{1-z/\rho}+O\left(1-z/\rho\right),\qquad
\lambda=\sqrt{2\rho+2\rho^2y'(\rho^2)}.
\end{equation}
Furthermore, the number~$y_n$ of trees of size~$n$ satisfies asymptotically
\begin{equation}\label{otter}y_n= \frac{\lambda}{2 \sqrt \pi} \cdot n^{-3/2}\rho^{-n} (1+O(1/n)).
\end{equation}
\end{lemma}
\begin{Proof}
Let~$C_{n-1}=\frac 1 n \binom{2n-2}{n-1}$ be the number of plane binary trees of size~$n$.
One has combinatorially $C_{n-1} 2^{-n}\le y_n\le C_n$. The~$C_n$,
which are the Catalan numbers and they are well-known to grow as $4^n n^{-3/2}$.
The bounds $1/4 \le \rho \le 1/2$ result.
It follows that $y(z^2)$ is analytic in a disc properly containing $|z|\le\rho$.
Then, from~\eqref{eq:gen_eq}, upon solving for~$y(z)$,
we obtain
\begin{equation}\label{eq:yz}
y(z)=1-\sqrt{1-2z-y(z^2)},
\end{equation}
which becomes singular when the argument of the square root vanishes.
The value $\rho$ is then the positive solution of
$2z+  y(z^2)=1$ and, at this point, we must have $y(\rho)=1$.
This reasoning also justifies the singular expansion~\eqref{singy},
seen to be valid in a $\Delta$-domain extending beyond the disc of convergence~$|z|<\rho$.

Equation~\eqref{otter} constitutes Otter's celebrated estimate:
it results from translating the square-root singularity of~$y(z)$ by means
of either Darboux's method or  singularity analysis~\cite{FlSe08,HaPa73,Otter48,Polya37}.
\end{Proof}

Numerically, one finds~\cite{Finch03,FlSe08,Otter48}:
\[
\rho \doteq 0.40269\,750367 
, \qquad
\lambda \doteq 1.13003\, 37163 
, \qquad
\frac{\lambda}{2\sqrt{\pi}}\doteq 0.31877\, 66259 
.
\]

\emph{Height.} Let $y_{h,n}$ be the number of trees of size~$n$ and
height \emph{at most}~$h$. 
(Height is measured as the maximum number
of edges along branches from the root to  external nodes.)
Let  $y_h(z)=\sum_{n\ge  1} y_{h,n}
z^n$. The arguments leading to (\ref{eq:gen_eq}) yield the fundamental recurrence
\begin{equation}\label{eq:gen_eq_h}
y_{h+1}(z)= z + \frac 1 2   y_h(z)^2 + \frac 1 2   y_h(z^2),
,\qquad \mbox{and}\qquad y_0(z)=z.
\end{equation}
We also set $
e_h(z)\equiv \sum_{n\ge 1} e_{h,n} z^n := y(z)-y_h(z),$ which is the generating function of 
trees with height \emph{exceeding} $h$. Then, a trite calculation shows that 
the $e_h(z)$ satisfy the main recurrence
\begin{equation}\label{eq:rec_eh}
e_{h+1}(z) = y(z) e_h(z) \pran{1-\frac{e_h(z)}{2y(z)}} + \frac 12 {e_h(z^2)},
\qquad \mbox{and} \qquad e_0(z)=y(z)-z,
\end{equation}
on which our subsequent treatment of height is entirely based.

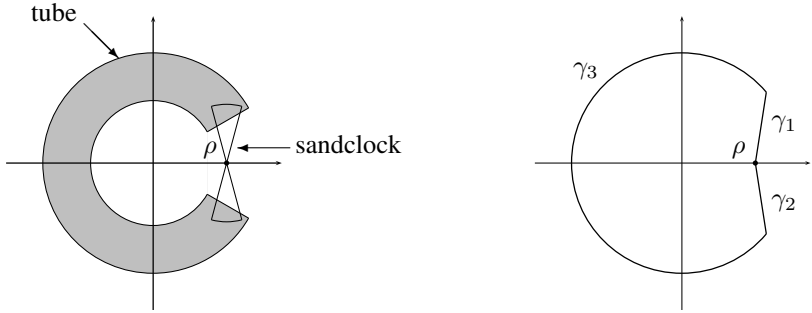
\begin{figure}[t]\centering
	\begin{picture}(300,120)
	\put(170,0){\includegraphics[width=4cm]{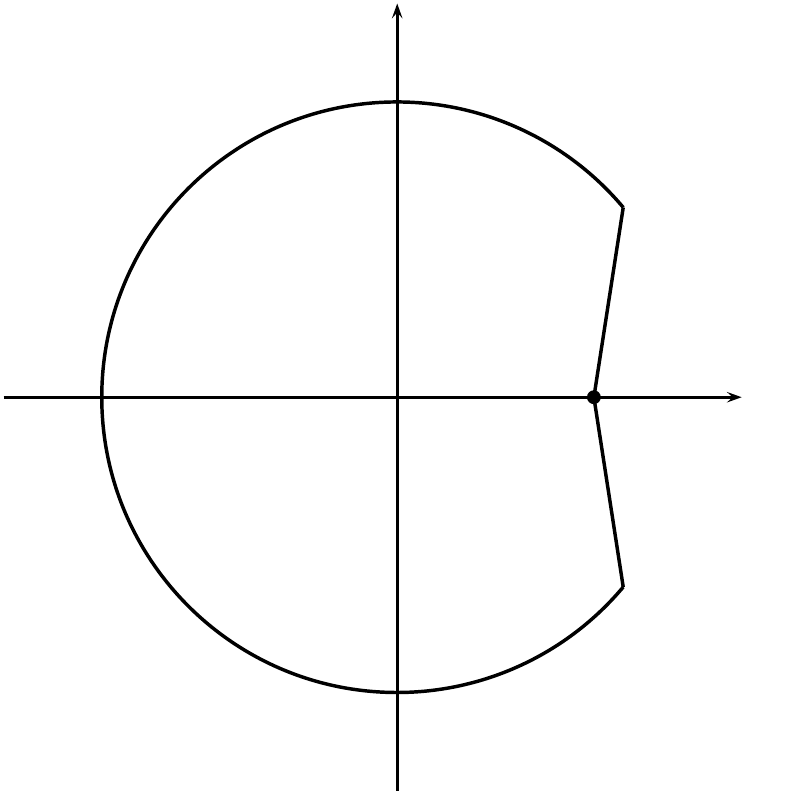}}
	\put(-30,0){\includegraphics[width=4cm]{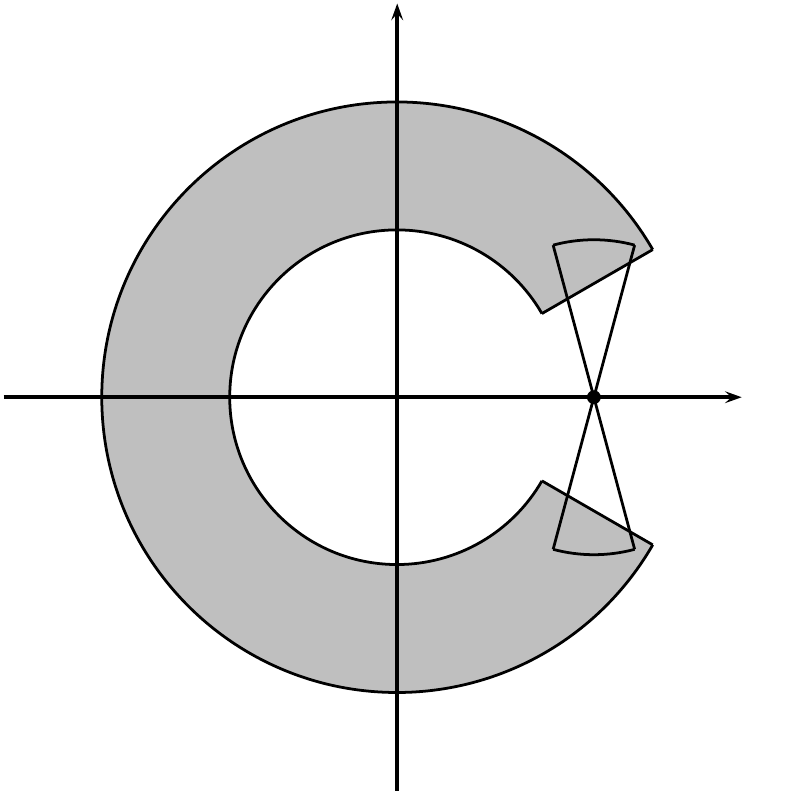}}
	\put(245,60){$\rho$}                       
	\put(45,60){$\rho$}
	\put(260,70){$\gamma_1$}              
	\put(260,40){$\gamma_2$}
	\put(185,90){$\gamma_3$} 
	\put(80,60){sandclock}   
	\put(78,62){\vector(-1,0){20}}       
	\put(-20,110){tube}   
	\put(00,109){\vector(1,-1){13}} 
        \end{picture}
\caption{\label{fig:hankel}The   ``tube''  and   ``sandclock''
	regions  on the left, and  the Hankel contour used to estimate
	$e_{h,n}$ on the right.}
\end{figure}

\emph{Analysis.} The distribution of height is accessible 
by
\begin{equation}\label{eq:gen_tail}
\p{H_n > h} = \frac{y_n - y_{h,n}}{y_n} = \frac{e_{h,n}}{y_n},
\end{equation} 
where $e_{h,n}=[z^n]e_h(z)$. 
We shall get a handle on its asymptotic properties by means of Cauchy's coefficient formula,
\begin{equation}\label{eq:cauchy}
e_{h,n} = \frac 1{2i\pi} \int_{\gamma} e_h(z)\frac{dz}{z^{n+1}},
\end{equation}
upon choosing a suitable integration contour~$\gamma$ in~\eqref{eq:cauchy}.
This task necessitates first developping suitable estimates of $e_h(z)$, for
values of~$z$ both \emph{inside} and \emph{outside} of the disc
of convergence~$|z|<\rho$. 
Precisely, we shall need estimates valid in a \emph{``tube''} around 
an arc of the circle~$|z|=\rho$,
\begin{equation}\label{eq:tube}
\mathcal T(\mu,\eta):=\{z:\quad -\mu <|z|-\rho<\mu ,~|\arg(z)|> \eta\}
.
\end{equation}
 as well as inside a \emph{``sandclock''}
anchored at~$\rho$ 
\begin{equation}\label{eq:sandclock}
\mathcal S(r_0,\theta_0):=\{z: \quad |z-\rho|< r_0,
\quad \pi/2-\theta_0 < |\arg(z-\rho)| < \pi/2+\theta_0\}.
\end{equation}
(See Figure~\ref{fig:hankel} for a rendering).

%

\emph{Plan.}
Estimates of  the  sequence of generating  functions~$(e_h(z))$ within
the disc  of  convergence and a tube,  where~$z$  stays away  from the
singularity~$\rho$ form the  subject of
Section~\ref{sec:away}.  The bulk of the technical work is relative to
the sandclock, in  Section~\ref{sec:singularity}.  We then develop our
main          approximation       in   Section~\ref{sec:estimates_eh}.
Sections~\ref{sec:away}--\ref{sec:estimates_eh}  closely follow    the
general strategy of  the original  paper~\cite{FlOd82}; however, nontrivial
adaptations are needed, due to the  presence of P\'olya terms, so that
the  problem is no  longer of a  ``pure'' iteration  type, as in \cite{FlOd82}.  We finally
reap  the   crop in Section~\ref{asy-sec},   where  our  main theorems
relative  to the   distribution of  the height  are  stated and proved
(these somewhat parallel  the local limit law  of~\cite{FlGaOdRi93} in
the planar case).

\section{Convergence away from the singularity}\label{sec:away}

Our aim in  this section\footnote{In what  follows, we freely omit the
arguments of $y$, $e_h$,  \dots, whenever they  are taken at $z$.}   is to
extend  the  domain where    $e_h$ is analytic     beyond the disc  of
convergence $|z|\le  \rho$,  when $z$  stays  in a ``tube''  $\mathcal
T(\mu,\eta)$ as defined in
\eqref{eq:tube} and is thus away from $\rho$. The main result is summarized by Proposition~\ref{prop:away},
at the  end of this  section.  Its proof  relies on the combination of
two   ingredients:    first,    the    fact,  expressed   by
Lemma~\ref{lem:bound1},   that $e_h\to 0$ (equivalently:
$y_h\to y$) in the  closed  disc of  radius~$\rho$; second,  a general
criterion for convergence  of  the $e_h$  to~0, which  is expressed by
Lemma~\ref{lem:criterion}.   The criterion implies  in essence that the
convergence domain is  an open set,  and  this fact provides  the basic
analytic continuation of the generating functions of interest.

\begin{lemma}\label{lem:bound1}
For all $z$ such that  $|z|\le \rho$, and $h\ge  1$, one has $|e_h(z)|\le \pran{|z|/\rho}^h/\sqrt h.$
\end{lemma}
\begin{proof}
To have height at  least $h$, a tree  needs  at least $h+1$  nodes, so
that $|e_h(z)| \le \sum_{n>h}y_n |z|^{n}.$ 
We first note an  easy  numerical refinement  of~\eqref{otter},  namely,  $y_n
\le\frac12 \rho^{-n} n^{-3/2}$. (See~\cite{FlGrKiPr95}  for  a detailed proof
strategy in the case of a similar but harder problem.) The claim follows by  bounding the sum using an integral.
\end{proof}

We now  devise a criterion for   the convergence of the $e_h$  to $0$.
This  criterion, adapted   from~\cite[Lemma~1]{FlOd82}, is  crucial in
obtaining   extended convergence   regions,    both near  the   circle
$|z|=\rho$ (in  this  section) and  near  the  singularity $\rho$  (in
Section~\ref{sec:singularity}).

\begin{lemma}[Convergence criterion]\label{lem:criterion}
Let $z\in \mathcal D:=\{z: |y|<1\}$ and assume that $|z|<\sqrt{\rho}$.
The  sequence $\{|e_h(z)|, h\ge  0\}$ converges  to 0  if and  only if
there exist an  integer  $m$ and real  numbers $\alpha,\beta\in(0,1)$,
such that the following three conditions are simultaneously met:
\begin{equation}\label{eq:criterion}
|e_m| < \alpha, \qquad |y|+\alpha /2< \beta, 
\qquad \alpha\beta+\pran{|z|^2/\rho}^m< \alpha.
\end{equation}
Furthermore, if \eqref{eq:criterion} holds then, for some $C$
and $\beta_0\in(0,1)$, one has $|e_h| \le C h \beta_0^h,$
for   all    $h\ge   m$.
\end{lemma}
\begin{Proof}
$(i)$ 
Since $|y|<1$ in $\mathcal D$ and $z<\sqrt \rho$, the convergence $e_h\to 0$ clearly implies that \eqref{eq:criterion} holds.

$(ii)$~Conversely, assume the three conditions in (\ref{eq:criterion}),
for some value~$m$. Then, they also hold for $m+1$. Indeed, recalling (\ref{eq:rec_eh}),
we see that, for any $h\ge 0$,
\begin{equation}\label{eq:bound_eh}
|e_{h+1}| ~\le~ |e_h| \pran{|y|+ \frac{|e_h|}2} + \frac{|e_h(z^2)|}2 ~\le~ |e_h| \pran{|y|+ \frac{|e_h|}2} + \pran{\frac {|z|^2}{\rho}}^h,
\end{equation}
where the P\'olya term  involving $|e_h(z^2)|$ has been bounded  using
Lemma~\ref{lem:bound1}.         Using     the  hypotheses
of~\eqref{eq:criterion}  together with \eqref{eq:bound_eh}  above taken at $h=m$,
yields $|e_{m+1}|<\alpha$.  
By induction, \eqref{eq:criterion} then holds for all $h\ge m$. 

$(iii)$~The assertion that $|e_h|\le Ch
\beta_0^h$ follows by expanding \eqref{eq:rec_eh} and using $\beta_0=\max\{\beta, |z|^2/\rho\}$.
\end{Proof}

We can now state the main convergence result of this section:

\begin{proposition}[Convergence in ``tubes'']\label{prop:away}
For  any  angle $\eta>0$, there  exists  a tube $\mathcal T(\mu,\eta)$
with width $\mu>0$,    such that $|e_h(z)| \to 0$,   as
$h\to\infty$, whenever $|z|$ lies in $\mathcal T(\mu,\eta)$.
\end{proposition}
\begin{proof}[sketch]
If we exclude a small sector of opening angle $2\eta$ around the
positive real axis, then the quantity, $
\lambda_0:=\sup \left\{\, |y(z)|;\quad |z|=\rho,~|\arg(z)|\ge \eta\, \right\},$
satisfies $\lambda_0<1$. 
%
The continuity of $y$ and Lemma~\ref{lem:criterion} then imply the convergence of $e_h(z)$ to $0$ in a small disc around each $z\in \{\rho e^{i\theta}: |\theta|\ge \eta\}$. The latter set being compact, we can extract a finite covering, which must then contain the desired tube where the convergence holds.
%
\end{proof}

\section{Convergence near the singularity}\label{sec:singularity}

We now focus on the behaviour of $e_h(z)$ in a ``sandclock'' around the singularity.
When $z$ approaches $\rho$, $|y|$ is no longer bounded away from 1 and the criterion for convergence given by Lemma~\ref{lem:criterion} cannot be used directly. However, as we prove below, the quantities $|e_h(z)|$ first exhibits a decreasing behaviour for $h\le N$, for some appropriate $N=N(z)$. At that point, $|e_N(z)|$ is small enough for the criterion of Lemma~\ref{lem:criterion} to be satisfied, whence enventually the convergence $|e_h(z)|\to 0$ as $h\to\infty$ in a sandclock.  

The upper bound we use for $|e_N|$ is based on the following alternative recurrence. 
\begin{lemma}[Alternative recurrence]\label{lem:rec_alt} If $e_i\ne 0$ and $e_i\left[1-e_i(z^2)/e_i^2\right]\ne 2y$ for $i=0,\dots, h-1$ then the following recurrence relations hold
\begin{equation}
\frac{y^h}{e_h} 
=\frac 1{e_0} + \frac 1 {2 y} \frac{1-y^h}{1-y} - \sum_{i=0}^{h-1} \frac{y^{i-1} e_i(z^2)}{2 e_i^2}
+ \sum_{i=1}^{h-1}  
\frac{y^{i-2} e_i}4\left[1-\frac{e_i(z^2)}{e_i^2}\right]^2
\pran{1-\frac{e_i}{2y}\left[1-\frac{e_i(z^2)}{e_i^2}\right]}^{-1}.\label{eq:alt_rec}
\end{equation}
\end{lemma}

\begin{proof}The proof relies on a classical idea in the study of slowly convergent iterations \cite{deBruijn81, FlOd82}.
Starting with the recurrence relation (\ref{eq:rec_eh}), rewritten as
$$\frac{e_{i+1}}{y^{i+1}} = \frac{e_i}{y^i} \pran{1-\frac{e_i}{2y} \left[1-\frac{e_i(z^2)}{e_i^2}\right]},$$
the trick is to \emph{take inverses}. This is a classical technique in the study of slowly convergent iterations (near an indifferent fixed point): see for instance \cite[page~152]{deBruijn81}.
Using the identity $(1-u)^{-1}=1+u+u^2(1-u)^{-1}$ shows that
$$
\frac{y^{i+1}}{e_{i+1}} - \frac{y^i}{e_i} = \frac {y^{i-1}}{2}\left[1-\frac{e_i(z^2)}{e_i^2}\right]  \pran{1-\frac{e_i}{2y}+ \frac{y^{i-2}}4 \left[1-\frac{e_i(z^2)}{e_i^2}\right]^2 \left[1-\frac{e_i(z^2)}{e_i^2}\right]}^{-1}. 
$$
Summing the terms of this equality for $i=0, \dots, h-1$ yields the claim. 
\end{proof}

An important step of the proof of Proposition~\ref{prop:bowtie} consists in controlling the behavior the terms in the alternative recurrence for all $h\le N$ (we will fix $N$ later).
In particular, we need good enough upper and lower bounds for $|e_h(z)|$, $h\le N$, and for $z$ around $\rho$ and $\rho^2$. 
Obtaining such estimates requires to study the recurrence relation (\ref{eq:rec_eh}) more carefully and to quantify the effect of the analytic P\'olya term $e_h(z^2)$ for $z$ in the ``sandclock''  $\mathcal S(r_0,\theta_0)$. The following lemma is evocative of the theory of iteration near an attractive fixed point and gives some estimates on the behaviour of $e_h$ in the interior of the disc of convergence. 


\begin{lemma}\label{lem:eh_rho2}There exist constants $K_1, K_2>0$ and $r>0$, such that, for all $h\ge 0$, and for $|z-\rho^2|<r$ one has $e_h(z) = C_h(z) \cdot y(z)^h,$
where $C_h(z)=(C(z)+o(1)) y(z)^h$ and $C(z)$ is analytic at $\rho^2$. Furthermore, $K_1 < |C_h(z)| < K_2$ and $\arg( e_h(z)) \le c_0 (h+1) r.$
\end{lemma}

With Lemma~\ref{lem:eh_rho2} in hand, we can obtain a first set of properties of $e_h(z)$ for $z\in \mathcal S(r_0,\theta_0)$ and $h$ not too large (depending on $z$). These will be used to derive an upper bound on $|e_N|$ and prove that $e_N$ satisfies the criterion of Lemma~\ref{lem:criterion}. In the following, we only need to consider $z\in \mathcal S(r_0,\theta_0)$, with $\Im(z)\ge 0$, since we clearly have $e_h(\bar z)=\overline{e_h(z)}$, where $\bar z$ denotes the complex conjugate of $z$.

\begin{lemma}[The initial behavior of $|e_h|$]\label{lem:init}Suppose $\Im(z)\ge 0$. Let $N(z)=\lfloor \arccos(1/4) / \arg y(z) \rfloor.$ There exist constants $r_0>0$ and $\theta_0>0$ such that if $z$ lies in the sandclock $ \mathcal S(r_0,\theta_0)$ then, for $1\le h\le N(z)$: 
\begin{equation}\label{eq:bound_arg}\frac{|y|^{h+1}}{2(h+1)}<|e_{h}(z)| < 1/ 2 \qquad \mbox{and}\qquad 0\le \arg(e_{h+1}) \le (h+1) \arg(y).\end{equation}
Furthermore, one has $|e_h(z)|<1/5$, for $6\le h \le N(z)$.
\end{lemma}
\begin{Proof}$(i)$ We first focus on the proof of the upper bound in \eqref{eq:bound_arg}. Consider the recurrence relation (\ref{eq:rec_eh}) rewritten as
\begin{equation}\label{eq:rec_ehy}
e_{h+1}/y = y \cdot e_h/y \pran{1-\frac {e_h/y} 2} + \frac{ e_h(z^2)}{2y}.
\end{equation}
The behavior of the first term in (\ref{eq:rec_ehy}) is dictated by properties of $g:w\mapsto w (1-w/2)$. A very similar function appeared in the analysis of \citet[Lemma 3]{FlOd82}. By a simple modification of the proof in \cite{FlOd82}, we have
\begin{equation}\label{eq:behav_g} \quad \left\{ \begin{array}{c} |w|\le 1 \\ 0 \le \arg w \le \arccos(1/4)\end{array}\right. \quad \Rightarrow \quad \left\{ \begin{array}{c} |g(w)|\le |w| \\ 0 \le \arg g(w) \le \arg w.\end{array}\right.
\end{equation}
We use \eqref{eq:behav_g} and an induction argument to prove that, for $0\le h\le N(z)$,
\begin{equation}\label{eq:induc_eh}|e_h|\le 1/2 \qquad \mbox{and} \qquad\arg e_h \le (h+1) \arg y.\end{equation} 

Write $z=\rho+re^{it}$. Then, 
by Lemma~\ref{lem:rho1},
we have  $e_0/y = 1 - z/y$ so that $|e_0/y| \le 1 - \rho + O(\sqrt r)$ and $0\le \arg(e_0/y) =O(\sqrt r).$
In particular, for $r$ small enough, \eqref{eq:induc_eh} holds for $h=0$.

Suppose now that \eqref{eq:induc_eh} holds for all integers up to $h$. 
To determine whether it also holds for $h+1$, we have to take the second term in the right-hand side of (\ref{eq:rec_ehy}) into account. 
For $r_0$ small enough and $z\in \mathcal S(r_0,\pi/8)$, Lemma~\ref{lem:eh_rho2} ensures that this second term cannot contribute any increase in the argument of $e_{h}/y$. 
Therefore, 
\begin{equation*}
\arg(e_{h+1}/ y) 
=\arg(e_h/y) + \arg (y) ~\le~(h+2) \arg(y).
\end{equation*}

Furthermore, since $y$ is analytic in $\mathcal D$, we have $|e_i| = |e_i(\rho)| +O(\sqrt r)$ as $r\to 0$, for all fixed $i\ge 0$. So for $r_0$ small enough and $z\in \mathcal S(r_0,\pi/8)$, $|e_{h+1}(z)|\le 1/2$ if $h\le 5$. Now, if $h\ge 5$, 
\[
|e_{h+1}/y| 
~\le~ |e_6(z)/y| + \frac 1 2 \cdot \sum_{i=6}^{h+1} |e_i(z^2)/y| \\
~\le~ |e_6(\rho)| + O(\sqrt r) + \frac 1 2 \cdot \sum_{i=6}^\infty (\rho + 3 r)^i,
\]
since $e_i(z^2)\le (|z|^2/\rho)^i$ by Lemma~\ref{lem:bound1}. One can then verify that, for $h\ge 6$, $
|e_{h+1}/y|
< 1/5,$
for $r_0$ small enough.
So, a fortiori, among the conditions in \eqref{eq:behav_g}, the one on the modulus holds as long as that on the argument does. 
The latter one holds for all $h\le N$.


$(ii)$ It remains to prove the lower bound in \eqref{eq:bound_arg}; we start with \eqref{eq:rec_ehy}.
For $h\le N(z)$, the additional P\'olya term $e_h(z^2)$ only contributes to making $|e_{h+1}|$ larger: for $z\in \mathcal S(r_0,\theta_0)$, by Lemma~\ref{lem:eh_rho2} and
the upper bound we have just proved, the arguments of both terms are such that, for $h<N(z)$,
$$|e_{h+1}/y| \ge |y| \cdot |e_h/y| \cdot \pran{1-\frac{|e_h/y|}2}.$$
Since $x\mapsto x (1-x/2)$ is increasing in $[0,1]$, we have, for all $h\ge 0$, $|e_h/y|\ge f_h$, where the sequence $(f_h)_{h\ge 0}$ is defined by  $f_{h+1}= |y| \cdot f_h \cdot \pran{1-f_h/2}$ and $f_0=|e_0|/|y|$.
The latter recurrence relation has been analysed by \citet{FlOd82} in the case of simply generated trees; the bound follows.
\end{Proof}

It is then easy to show that, as expected, the P\'olya terms are essentially negligible:

\begin{lemma}\label{lem:polya_term}There exist $r_0>0$ and $\theta_0>0$ small enough such that, for $z\in \mathcal S(r_0,\theta_0)$ and for all $h\le N$, one has 
$|e_h(z^2)|/ |e_h(z)^2|\le  \min\left\{12  \cdot (2\rho)^{h+1}, \frac 1 2 \right\}.$
\end{lemma}

Finally, we prove the main result of this section. The proof follows the lines of the analogous statement \cite[Proposition 4]{FlOd82}, where the iteration is ``pure'', but now needs to control the effect of the P\'olya terms, using Lemma~\ref{lem:polya_term}. 
\begin{proposition}[Convergence in a ``sandclock'' around $\rho$]\label{prop:bowtie}There exist constants $r_0$ and $\theta_0$ such that the sequence $\{e_h(z), h\ge 0\}$ converges to zero for $z\in\mathcal S(r_0,\theta_0)$, where the ``sandclock'' $\mathcal S(r_0,\theta_0)$ is defined by \eqref{eq:sandclock}.
\end{proposition}
\begin{proof} We prove that for $h=N\equiv N(z)$ defined in Lemma~\ref{lem:init}, $e_N$ satisfies the convergence criterion of Lemma~\ref{lem:criterion}. For this purpose, we use the alternative recurrence stated in Lemma~\ref{lem:rec_alt}
\begin{equation}\label{eq:rec_labeled}
	\frac{y^h}{e_h} 
	= \frac 1 {2 y} \frac{1-y^h}{1-y} + \underbrace{\frac 1{e_0}- \sum_{i=0}^{h-1}\frac{y^{i-1} e_i(z^2)}{2 e_i^2}}_{A_h} 
	+ \underbrace{\sum_{i=1}^{h-1}  
	{\frac{y^{i-2} e_i}4\left[1-\frac{e_i(z^2)}{e_i^2}\right]^2}
	\pran{1-\frac{e_i}{2y}\left[1-\frac{e_i(z^2)}{e_i^2}\right]}^{-1}}_{B_h}
\end{equation}
and devise an asymptotic lower bound for the right-hand side.
First observe that we can indeed use the relation, since by Lemmas~\ref{lem:init} and \ref{lem:polya_term}, for all $i=0, \dots, N$, the denominators do not vanish.

Write $1-y(z)=\epsilon e^{it}$. Obtaining the claim reduces to proving properties for small $\epsilon>0$ and $t$ close to $-\pi/4$. 
The following expansions are valid uniformly for $t\in [-\pi/4 -\delta, -\pi/4 +\delta]$ with $0<\delta<\pi/4$ when $\epsilon \to0$:
\begin{equation}\label{eq:value_N}
\left\{\begin{array}{r l l}
1-|y|&=&\epsilon \cos t +O(\epsilon^2),\\
\arg (y) &=& -\epsilon \sin t +O(\epsilon^2),
\end{array}\right.
\qquad \mbox{and} \qquad
\left\{\begin{array}{rl l}
	N(z) &=& -\varphi/(\epsilon \sin t) + O(1)\\
	1-|y|^N &=& 1- e^{\varphi \cdot \cot t} + O(\epsilon),
\end{array}\right.
\end{equation}
where $\varphi:=\arccos(1/4)$.
The first term of the right-hand side of \eqref{eq:rec_labeled} brings the main contribution:
$$\left|\frac 1 {2 y} \frac{1-y^N}{1-y}\right|= \frac 1 {2|y|}\cdot \frac{|1-y^N|}{|1-y|}\ge \frac1 2 \frac{1-|y|^N}{|1-y|}= \frac{1-e^{\varphi \cdot \cot t}}{2\epsilon}+O(1),$$
as $\epsilon\to0$. On the other hand, we have $|A_N|=O(1)$
for, by Lemma~\ref{lem:polya_term}, the summands decrease geometrically.
The second error term $B_N$ appearing in (\ref{eq:rec_labeled}) can be bounded by splitting the sum into two at $h=K$. Let $\nu>0$. By Lemma~\ref{lem:polya_term}, there is $K\ge 6$ large enough that for all $i$ satisfying $K\le i \le N$, $|e_i(z^2)|/|e_i(z)^2| \le \nu.$
Then, by Lemma~\ref{lem:init}, and for $\epsilon$ small enough, we have
\begin{equation}\label{eq:bN}
|B_N| \le \frac {K/2} 4 \cdot \frac{(3/2)^2}{1-\frac{1/2}{2-2\epsilon} \cdot (3/2)} + \frac{1/5} 4 \cdot \frac{(1+\nu)^2}{1-\frac{1/5}{2-2\epsilon} (1+\nu)} \frac{1-|y|^N}{1-|y|}< \frac K 2 + \frac 3 {50} \cdot \frac{1-|y|^N}{1-|y|},
\end{equation}
upon choosing $\nu=1/100$.
It follows that
\begin{equation}\label{eq:ynen}\frac{|y^N|}{|e_N|} \ge \frac{1-e^{\varphi \cdot \cot t}}\epsilon \pran{\frac 12 -\frac 3 {50 \cos t}} + O(1) > \frac 2 5 \cdot \frac{1-e^{\varphi \cdot \cot t}}\epsilon,
\end{equation}
for all $z\in \mathcal D$ such that $\epsilon<\epsilon_0$ and $|t-\pi/4|< \delta_0$, as soon as both $\epsilon_0$ and $\delta_0$ are small enough.

We can now focus on the criterion for convergence (Lemma~\ref{lem:criterion}) and the conditions in \eqref{eq:criterion} for $m=N$. From \eqref{eq:value_N} and (\ref{eq:ynen}) we have
for $\epsilon>0$ small enough, 
\begin{equation*}\label{eq:crit2}
|y| + \frac{|e_N|} 2 \le 1- \epsilon \cdot \pran{\cos t  - \frac{5}{4} \cdot \frac {e^{\varphi \cdot \cot t}} {1-e^{\varphi \cdot \cot t}}} +O(\epsilon^2).
\end{equation*}
A simple verification shows that the coefficient of $\epsilon$ above is at most $-1/4$ for all $t$ close enough to $-\pi/4$.
Thus, for all $\epsilon>0$ small enough,
\begin{equation}\label{eq:crit_values}|e_N| \le\frac{5}{2} \cdot \frac{\epsilon \cdot  e^{\varphi \cot t}}{1-e^{\varphi \cdot \cot t}}=: \alpha \qquad \mbox{and} \qquad |y| + \frac{|e_N|}2 < 1- \frac {\epsilon}{4}=: \beta.\end{equation}
Also, by (\ref{eq:value_N}), we have $(|z|^2/\rho)^N=o(\alpha(1-\beta))$. 
Hence the criterion for convergence of Lemma~\ref{lem:criterion} is satisfied with the values of $\alpha$ and $\beta$ specified in (\ref{eq:crit_values}), as soon as $\epsilon$ is small enough. 
\end{proof}

\section{Main approximation}\label{sec:estimates_eh}

In this section, we estimate $e_h(z)$ around the singularity.

\begin{proposition}[Main estimate for $e_h$ in a sandclock around $\rho$]\label{prop:eh_sim}There exist $r_1, \theta_1$ and $K, K'$ such that for all $z\in \mathcal S(r_1,\theta_1)$ and all $h\ge 1$,
$$\frac{y^h}{e_h} = \frac 1 2 \cdot \frac{1-y^h}{1-y} + R_h(z) \qquad \mbox{where} \qquad |R_h(z)| \le K \min\left\{\log \frac 1{1-|y|}, \log(1+h)\right\}.$$
Furthermore, $|R_h-R_{h+1}|\le K'/h$.
\end{proposition}

The main idea is to obtain a better control on the error term using bounds extending those obtained in Section~\ref{sec:singularity} for $h>N$, knowing the \emph{a priori} information that $e_h$ converges. The proof of Proposition~\ref{prop:eh_sim} also requires the bounds to be \emph{uniform} in the distance to the singularity $|z-\rho|$ and in $h$. 
\begin{lemma}[Uniform lower bound for $|e_h|$]\label{lem:low2}For any $\delta\in(0,1)$, there exist constants $r_1,\theta_1>0$ such that if $z\in \mathcal S(r_1,\theta_1)$, then for all $h\ge 0$, one has 
$|e_h(z)|\ge (1-\delta)^{h+2}/(2(h+1)).$
\end{lemma}
\begin{Proof}Let $\delta\in(0,1)$. For $r$ small enough, we have $|y|>1-\delta/2$ provided $r:=|z-\rho|<r_0$. Then, by Lemma~\ref{lem:init}, the desired lower bound is clear for $h\le N$. So we now assume that $h> N$. The recurrence relation (\ref{eq:rec_eh}) implies that
$$|e_{h+1}|\ge |y| |e_h| \pran{1-\frac{|e_h|}{2|y|}} - \frac{|e_h(z^2)|}2.$$
By Lemma~\ref{lem:criterion}, $|e_h|$ is decreasing for $h\ge N$, when the criterion is satisfied. So, for $z\in \mathcal S(r_0,\theta_0)$ as in Proposition~\ref{prop:bowtie}, we have
$$
|e_{h+1}|\ge |y| \pran{1-\frac{|e_N|}{2|y|}} \cdot |e_h| - \frac{|e_h(z^2)|}2.$$
However, we have set $N$ in such a way that $|y|+ |e_N|/2<1$, and by Lemma~\ref{lem:bound1}, we have
$|e_{h+1}|\ge (1-\delta) |e_h| - (\rho+r)^h.$
Routine verifications then show, for $r$ small enough, that $1-\delta>\rho+r$, the negative contribution decreases fast enough so that $|e_h|$ remains bounded from below as desired.
\end{Proof}

We can now proceed with the proof of a uniform upper bound for $|e_h|$ when $z\in \mathcal S(r_1,\theta_1)$.

\begin{lemma}[Uniform upper bound for $|e_h(z)|$ around $\rho$]\label{lem:up_uniform}There exist constants $c_1$, $r_1$ and $\theta_1$ such that, for any $h\ge 1$, and $z\in \mathcal S(r_1,\theta_1)$, we have $|e_h(z)|\le c_1/h.$
\end{lemma}
\begin{Proof}
Write $1-y=\epsilon e^{it}$ for some $\epsilon>0$ and $t$. It suffices to prove that the result holds for all such $z$ provided $\epsilon$ is small enough and $t$ close enough to $-\pi/4$. We use again \eqref{eq:rec_labeled}. Then, we have, 
\begin{equation}\label{eq:unif1}
\left|\frac 1 {2 y} \frac{1-y^h}{1-y}\right|= \frac 1 {2|y|}\cdot \frac{|1-y^h|}{|1-y|}\ge \frac{1-|y|^h}{2|1-y|}= \frac{1-|y|^h}{2\epsilon}.
\end{equation}
The error terms $A_h$ and $B_h$ are bounded as in the proof of Proposition~\ref{prop:bowtie} and we obtain, for all $h\ge 0$ and $\epsilon>0$ small enough,
$|A_h|\le K_1$ and $|B_h|\le K_2 + \frac 1 5(1-|y|^h)/(1-|y|)$\footnote{In what follows, we use generically~$K,K_1,\ldots$ to denote 
positive absolute constants, not necessarily of the same value at different occurrences.}.
One then sees that, for all $h\ge 0$, and all $t$ close enough to $-\pi/4$,
\begin{equation}\label{eq:eh_unif}
\frac{|y^h|}{|e_h|} \ge \frac{1-|y|^h}{\epsilon} \pran{\frac 1 2 - \frac 1{5 \cos t}} - K_3> \frac{1-|y|^h}{5\epsilon}  - K_3 \qquad \mbox{and}\qquad |e_h| \le 10\epsilon \frac{|y|^h}{1-|y|^h}. 
\end{equation}
For $h$ not too large, $|y|^h$ decreases at least linearly in $h$, and one can show that $|y|^h\le 1-\delta h \epsilon$, for some small $\delta>0$ and all $0\le h\le N$ as long as $|t+\pi/4|$ is small enough. Equation \eqref{eq:eh_unif} then implies that $|e_h|\le K_4/h$, for all $h\le N$.
 
On the other hand, by (\ref{eq:value_N}), if $h\ge N$ then the factor $|y|^h$,  ensures the desired decreasing behaviour. Indeed, for $\epsilon$ small enough and $t$ close enough to $-\pi/4$, one has $|y|<(1-\epsilon/2)$ and $|e_h| \le K_5 \epsilon (1-\epsilon/2)^h.$ The maximum of the right-hand side above is obtained for $\epsilon=2/(h+1)$, which implies that
$|e_h| \le 2K_5/(h+1)$ for $h \ge N.$ This completes the proof.
\end{Proof}

\begin{proof}[of Proposition~\ref{prop:eh_sim}]The proof consists in using Lemma~\ref{lem:up_uniform} above to bound the error terms in \eqref{eq:rec_labeled} for $z\in \mathcal S(r_2,\theta_2)$, with $r_2=\min\{r_0,r_1\}$ and $\theta_2=\min\{\theta_0,\theta_1\}$. 
For some constants $c_2$ and $c_3$, we have
\[
|A_h|+|B_h|\le\frac {11}{(1-2\rho)^3}+c_2 \pran{1 + \sum_{i=1}^{h}\frac{|y^i|}{i}} \le  c_3\min\left\{\log\pran{\frac 1{1-|y|}}, 1+\log h \right\}.
\]
Also, since $A_h$ and $B_h$ are partial sums, 
$R_h-R_{h+1}$ only contains one summand, 
which is easily seen to be $O(1/h)$ uniformly by Lemmas~\ref{lem:low2} and \ref{lem:up_uniform}.
\end{proof}

\section{Asymptotic analysis and distribution estimates}\label{asy-sec}

The basis of our estimates relative to the distribution of height is 
Proposition~\ref{prop:eh_sim}
in conjunction with Cauchy's coefficient formula~\eqref{eq:cauchy}
where~$\gamma$ is a contour 
comprised of the arc of an outer circle of radius larger than~$\rho$ 
(and interior to the cardioid-shaped 
region, where $|y|<1$) and a set of two connecting segments
passing through the singularity~$\rho$  (Figure~\ref{fig:hankel}).
 In addition, it proves useful, in order to
garantee well-defined determinations of square roots, to think of the two segments as in fact joined by 
an infinitesimal arc of a circle that passes to the \emph{left} of the singularity~$\rho$.
The strategy just described belongs to the general orbit of singularity analysis methods.

\begin{theorem}[Limit law of height] \label{clt}
The height $H_n$ of a random tree taken uniformly from $\cal Y_n$ admits
a limiting theta distribution: for any fixed~$x>0$, there holds
\[
\lim_{n\to\infty} \Pr(H_n\ge x\sqrt{n})=\sum_{k\ge1} 
(k^2\lambda^2 x^2-2)e^{-k^2\lambda^2 x^2/4},
\qquad \lambda:=\sqrt{2\rho+2\rho^2y'(\rho^2)}.
\]
\end{theorem}
\begin{proof}
The integration contour in Cauchy's formula~\eqref{eq:cauchy}
is $\gamma=\gamma_1\cup\gamma_2\cup\gamma_3$. There, $\gamma_1$ is the segment lying in
the half-plane $\Im(z)\ge0$, $\gamma_2$ is the complex-conjugate image of~$\gamma_1$, and
$\gamma_3$ is the outer circular arc. For the radius of the latter circle, we adopt
$r_n=\rho(1+\log^2 n/n).$
We assume that $\gamma_3$ lies in a legal tube (granted by Proposition~\ref{prop:away}) and that $\gamma_1$ and $\gamma_2$ are in an overlapping sandclock such that Proposition~\ref{prop:eh_sim} applies. We set, for some $\theta_1>0$ :
$\gamma_1=\bar\gamma_2=\left\{\rho\pran{1+x e^{i\pi/2+i\theta_1}}: x\in [0,\delta_n]\right\},$ with $\delta_n$ such that the rectilinear portions $\gamma_1$ and $\gamma_2$ connect with the outer circle $\gamma_3$. So, we have 
$\delta_n \sim \log^2 n/(n \sin \theta_1).$

\emph{Outer circular arc $(\gamma_3)$.}
By Proposition~\ref{prop:away}, we have $e_h(z)\to0$ uniformly on $\gamma_3$ as $h\to\infty$. In particular, all moduli $|e_h(z)|$
are bounded by an absolute
 constant $K$. On the other hand the Cauchy kernel $z^{-n}$ is small
on the outer part of the contour, so that
\begin{equation}\label{outer}
\left|\int_{\gamma_3} e_h(z)\, \frac{dz}{z^{n+1}}\right|<K_1\rho^{-n} \exp\left(-\log^2 n\right).
\end{equation}
This contribution is thus exponentially small compared to~$y_n$.

\emph{Rectilinear parts $(\gamma_1,\gamma_2)$.} Our 
objective is to replace~$e_h$ by the simpler quantity
\begin{equation}\label{wheh}
\wh e_h(z) \equiv \wh e_h := 2\frac{1-y}{1-y^h} y^h,
\end{equation}
as suggested by Proposition~\ref{prop:eh_sim}.
Along~$\gamma_1,\gamma_2$, the singular expansion of $y(z)$ applies, so that
$1-y=O((\log n)/n^{1/2})$ and the error term~$R_h(z)$ is $O(\log n)$.
There results that $(1-y^h)/(1-y)$ is always 
 at least as large in modulus as $K_2\sqrt{n}/\log n$ (by a study of the variation of
$|1-e^{-h\tau}|/|1-e^{-\tau}|$), and we have
\begin{equation}\label{ap1}
\frac{y^h}{e_h}= 
\frac{y^h}{\wh e_h}\left(1+O\left(\frac{\log^2n}{\sqrt{n}}\right)\right).
\end{equation}

It proves convenient  to define the following approximation for $e_{h,n}$:
\begin{equation}\label{Ehn}
E(h,n):=\frac{1}{2i\pi} \int_{\gamma_1\cup\gamma_2} \wh e_h \, \frac{dz}{z^{n+1}},
\end{equation}
and effect the change of variables
$z=\rho\left(1-\frac{t}{n}\right)$.
The quantity~$t$ then varies from $-i n \delta_n e^{-i\theta_1}$, loops to the right of the origin, then moves away  to $i n\delta_n e^{i\theta_1}$.
With the singular expansion of~$y(z)$ as in~\eqref{singy},
we have on~$\gamma_1,\gamma_2$,
\begin{equation}\label{ap11}
z^{-n}=\rho^{-n} e^t \left(1+O\left(\log^4 n/ n\right)\right),
\qquad
y(z)= 1- \lambda \sqrt{t/n}+O\left(t/n\right).
\end{equation}
and, with $h=x\sqrt{n}$ and $|t|\le K_2 \log^2n$:
\begin{equation}\label{ap12}
y^h 
= \exp(-\lambda x\sqrt{t}) \cdot \left(1+O\left(\log^2 n/\sqrt{n}\right)\right).
\end{equation}
We also find\footnote{The expression $\log^\star n$ represents an unspecified positive
power of $\log n$.}, for the range of values of~$t$ corresponding to~$\gamma_1,\gamma_2$: 
\begin{equation}\label{ap2}
\frac{1-y^h}{1-y}
=\left[\sqrt{n}\cdot \frac{1-\exp(-\lambda x\sqrt{t})}{\lambda\sqrt{t}}\right]
\left(1+O\left(\frac{\log^\star n}{\sqrt{n}}\right)\right).
\end{equation}

The approximations~\eqref{ap11},~\eqref{ap12}, and~\eqref{ap2}
motivate considering, as an approximation of~$E(h,n)$ in~\eqref{Ehn}, the contour integral
\begin{equation}\label{Ix}
J(X):=\cau\int_\cal L \frac{\exp(-X\sqrt{t})}{1-\exp(-X\sqrt{t})}
\sqrt{t} e^t\, dt =\cau \sum_{k\ge1}\int_\cal L {\exp(-kX\sqrt{t})} \sqrt{t}
e^t\, dt,
\end{equation}
where $\cal L$ is the image of $\gamma_1\cup \gamma_2$ in the change of variable.
We can now make $J(X)$ explicit.
Each integral on the right side can be evaluated by the change of variables $w=i\sqrt{t}$, $t=-w^2$. 
By completing the square and flattening the image contour~$\cal L$ onto the real line,
we obtain:
\begin{equation}\label{Jx}
J(X)=\frac{1}{4\sqrt{\pi}}\sum_{k\ge1} e^{-k^2X^2/4}(k^2X^2-2).
\end{equation}

\emph{Error management.} 
%
It can be checked that the replacements: $e_h\mapsto \wh e_h$, $y\mapsto 1-\lambda \sqrt{t/n}$, and $y^h\mapsto\exp(-\lambda x \sqrt t)$
only entail error terms of order $\rho^{-n} n^{-2}\log^\star n$, which implies, for $h=x\sqrt n$,
%
for $h=x\sqrt{n}$:
\[
e_{h,n}=2\lambda \rho^{-n}n^{-3/2} J(\lambda x)
+O\left(\rho^{-n} \log^\star n/n^2\right).
\]
The explicit form of~$J(X)$ in~\eqref{Jx} and
the asymptotic form of~$y_n$ (Lemma~\ref{lem:rho1}),
yield
the statement.
\end{proof}

\begin{theorem}[Local limit law of height] \label{llt}
The height $H_n$ of a random tree taken uniformly from $\cal Y_n$ admits
a local limiting distribution: for $x$ in a compact set of $\R_{>0}$ 
and $h=x\sqrt{n}$
an integer, there holds
\[
\Pr(H_n=h)\sim\frac{1}{2x\sqrt{n}}
\sum_{k\ge1} 
(k^4\lambda^4 x^4-6k^2\lambda^2 x^2)e^{-k^2\lambda^2 x^2/4}.
\]
\end{theorem}
\begin{Proof}
Proceeding like in the proof of Theorem~\ref{clt},
we can justify approximating the number of trees
of height exactly~$h$ by the integral
\[
\frac{1}{2i\pi}\int_\gamma \left(\wh e_h-\wh e_{h+1}\right)\, \frac{dz}{z^{n+1}},
\qquad \mbox{where} \qquad
\wh e_h - \wh e_{h+1}= 2y^h \frac{(1-y)^2}{(1-y^h)(1-y^{h+1})}.
\]
The approximations~\eqref{ap12} and~\eqref{ap2} then motivate consideration of the
integral
\[
J_1(X):=\cau\int_\cal L \frac{\exp(-X\sqrt{t})}{(1-\exp(-X\sqrt{t}))^2}
t e^t\, dt,
\]
and one finds (with the auxiliary estimate~$R_h-R_{h+1}=O((\log^\star n)/\sqrt{n})$ provided by Proposition~\ref{prop:eh_sim}):
\[
y_{n,h}-y_{n,h+1}=2\lambda^2 \rho^{-n} n^{-2} 
J_1(\lambda x)+O\left(\rho^{-n}\frac{\log^\star n}{n^{5/2}}\right).
\]
On the other hand, differentiation under the integral sign yields $J_1(X)=-J'(X)$,
which proves the statement.
\end{Proof}

Revisiting the proof of Theorems~\ref{clt} and~\ref{llt}
shows that one can allow $x$ to become either small or large,
albeit to a limited extent. Indeed, it can be checked,
for instance, that allowing $x$ to get as large as $O(\sqrt{\log n})$
only introduces extra powers of~$\log n$ in error estimates. 
However, such extensions are limited by the fact that the main theta term eventually becomes
smaller than the error term. We state (compare with~\cite[Th.~1.1]{FlGaOdRi93}):

\begin{theorem}[Moderate deviations]\label{thm:moderate}
There exist constants $A,B,C>0$ such that for $h=x\sqrt{n}$ with
$A/\sqrt{\log n}\le x\le A\sqrt{\log n}$ and $n$ large enough, there holds
\begin{equation}\label{mld}
\left|\Pr(H_n\ge x\sqrt{n})-\sum_{k\ge1}(k^2x^2-2)e^{-k^2x^2/4}
\right|
\le \frac{C}{n^B}.
\end{equation}
In particular, if~$x\to\infty$ in such a way that  $x\le A\sqrt{\log n}$,
then, uniformly,
\[
\Pr(H_n\ge x\sqrt{n})
\sim \lambda^2 x^2 e^{-\lambda^2x^2/4}.
\]
Similar estimates hold for the local law.
\end{theorem}

These estimates can be supplemented by (very) large deviation estimates in the style of~\cite[Th.1.4]{FlGaOdRi93}: 
it suffices to make use of the fact that~$e_h$ 
is bounded from above by a large power and optimize on~$r\in(0,\rho)$ the
saddle-point bound
\[
e_{h,n} \le \frac{e_h(r)}{r^n},\qquad 0<r\le \rho.
\]
The probability of a linearly height is then exponentially small:

\begin{theorem}[Very large deviations]\label{thm:large_deviations}
There exists a continuous decreasing function~$I(\xi)$ satisfying $I(\xi)>0$ for $0<\xi\le 1$,
such that, given any fixed~$\delta>0$, 
one has for $h=\xi n$, and $\delta<\xi<1-\delta$, 
\[
\Pr(H_n\ge \xi n)\le K n^{3/2} e^{-n I(\xi)},
\]
where~$K$ depends on~$\delta$.
\end{theorem}
Finally, the approximation of~$e_h$ by~$\wh e_h$ in~\eqref{wheh} is
good enough to grant us access to moments (cf also~\cite{FlOd82}). The problem reduces to estimating generating functions
of the form 
\[
M_r(z)=2 (1-y)^2 \sum_{h\ge1}h^r \frac{y^h}{(1-y^h)^2},
\]
which are accessible to the Mellin transform technology~\cite{FlGoDu95}, upon
setting~$y=e^{-\tau}$.

\begin{theorem}[Moments of height]\label{thm:moments}
Let~$r\ge1$. The $r$\emph{th} moment of the height~$H_n$ satisfies
\[
\E{H_n}\sim \frac 2 \lambda \sqrt{\pi n}
\qquad\mbox{and}\qquad
\mathbb{E}[H_n^r]\sim r(r-1)\zeta(r) \Gamma(r/2) \left(\frac{2}{\lambda}\right)^r n^{r/2},\quad r\ge 2
.\]
\end{theorem}


\begin{small}
\setlength{\bibsep}{.3em}
\bibliographystyle{abbrvnat}
\bibliography{algo}
\end{small}

\end{document}